\documentstyle[twoside,epsfig]{article}
\input epsf
\title{A note on a series containing the Laguerre polynomials}
\author{\sc Y. S. Kim$^a$,  A. K. Rathie$^b$ and R. B. Paris$^c$\\
\\
${}^a\!$ Department of Mathematics Education, Wonkwang University, Iksan, Korea\\
E-Mail: yspkim@wonkwang.ac.kr\\
${}^b\!$ Department of Mathematics, Central University of Kerala, Kasaragad 671123,\\
Kerala, India\\
E-Mail: akrathie@cukerala.edu.in\\
${}^c\!$ School of Engineering, Computing and Applied Mathematics,\\
 University of Abertay Dundee, Dundee DD1 1HG, UK\\
E-Mail: r.paris@abertay.ac.uk}
\begin{document}
\def\f#1#2{\mbox{${\textstyle \frac{#1}{#2}}$}}
\def\dfrac#1#2{\displaystyle{\frac{#1}{#2}}}
\def\boldal{\mbox{\boldmath $\alpha$}}
\newcommand{\bee}{\begin{equation}}
\newcommand{\ee}{\end{equation}}
\newcommand{\lam}{\lambda}
\newcommand{\ka}{\kappa}
\newcommand{\al}{\alpha}
\newcommand{\th}{\theta}
\newcommand{\om}{\omega}
\newcommand{\Om}{\Omega}
\newcommand{\fr}{\frac{1}{2}}
\newcommand{\fs}{\f{1}{2}}
\newcommand{\g}{\Gamma}
\newcommand{\br}{\biggr}
\newcommand{\bl}{\biggl}
\newcommand{\ra}{\rightarrow}
\newcommand{\mbint}{\frac{1}{2\pi i}\int_{c-\infty i}^{c+\infty i}}
\newcommand{\mbcint}{\frac{1}{2\pi i}\int_C}
\newcommand{\mboint}{\frac{1}{2\pi i}\int_{-\infty i}^{\infty i}}
\newcommand{\gtwid}{\raisebox{-.8ex}{\mbox{$\stackrel{\textstyle >}{\sim}$}}}
\newcommand{\ltwid}{\raisebox{-.8ex}{\mbox{$\stackrel{\textstyle <}{\sim}$}}}
\renewcommand{\topfraction}{0.9}
\renewcommand{\bottomfraction}{0.9}
\renewcommand{\textfraction}{0.05}
\newcommand{\mcol}{\multicolumn}
\date{}
\maketitle
\begin{abstract}
Expressions for the summation of a new series involving the Laguerre polynomials
are obtained in terms of generalized hypergeometric functions.
These results provide alternative, and in some cases simpler, expressions to those recently obtained in the literature.
\vspace{0.4cm}

\noindent {\bf Mathematics Subject Classification:} 33C15, 33C20 \setcounter{section}{1}
\setcounter{equation}{0}
\vspace{0.3cm}

\noindent {\bf Keywords:} Laguerre polynomials, generalized hypergeometric functions, generalized Kummer summation theorem 
\end{abstract}
\vspace{0.3cm}

\begin{center}
{\bf 1. \  INTRODUCTION}
\end{center}
\setcounter{section}{1}
\setcounter{equation}{0}
\renewcommand{\theequation}{\arabic{section}.\arabic{equation}}
For arbitrary real $\nu$ the polynomial solutions of the differential equation
\[xy''+(\nu+1-x)y'+ny=0\]
are called the generalized Laguerre polynomials and are usually denoted by $L_n^{(\nu)}(x)$. The Laguerre polynomials are encountered in many branches of pure and applied mathematics and mathematical physics, especially in quantum mechanics and in the radial part of the solution of the Schr\"odinger equation for a single-electron atom. They form an orthogonal set on $[0, \infty)$ with the weight function $x^\nu e^{-x}$, with the first three polynomials given by 
\[L_0^{(\nu)}(x)=1, \quad L_1^{(\nu)}(x)=1-x+\nu, \quad  L_2^{(\nu)}(x)=\fs x^2-(\nu+2)x+\fs(\nu+1)(\nu+2).\]
In general, $L_n^{(\nu)}(x)$ can be represented as a terminating confluent hypergeometric function ${}_1F_1$ in the form
\[L_n^{(\nu)}(x)=\frac{(\nu+1)_n}{n!}\,{}_1F_1(-n; \nu+1; x).\]
Here $(a)_n$ denotes the Pochhammer symbol, or rising factorial, defined by $(a)_n=\g(a+n)/\g(a)$.

In \cite{KRP}, Kim {\it et al.\/} obtained summation formulas for the series involving the generalized Laguerre polynomial $L_n^{(\nu)}(x)$ given by
\[\sum_{n=0}^\infty \frac{x^n\,L_n^{(\nu)}(x)}{(1\pm\nu+j)_n}\]
for integer $j$, where $-5\leq j\leq 5$, following on from an earlier investigation by Exton \cite{E} in the case $j=0$. Recently, Brychkov \cite{B} has extended these results for any integer $j$. The aim of this note is to derive alternative expressions for the summation of the series
\[S(\pm\nu, \pm j)\equiv e^{-x}\sum_{n=0}^\infty \frac{x^n\,L_n^{(\nu)}(x)}{(1\pm\nu\pm j)_n}\]
for any non-negative integer $j$. Our results are different from, and in some cases simpler than, those obtained in \cite{B} and \cite{KRP}.
\vspace{0.6cm}

\begin{center}
{\bf 2. \  THE SERIES $S(\nu, \pm j)$}
\end{center}
\setcounter{section}{2}
\setcounter{equation}{0}
\renewcommand{\theequation}{\arabic{section}.\arabic{equation}}
We start with the transformation \cite[(3.5)]{KRP}
\bee\label{e1}
e^{-x} \sum_{n=0}^\infty \frac{(a_1)_n \ldots (a_p)_n}{(b_1)_n \ldots (b_q)_n}\,(-xy)^n\,L_n^{(\nu)}(x)=\sum_{n=0}^\infty \frac{(-x)^n}{n!}\,{}_{p+2}F_q\left[\!\!\begin{array}{c}-n, -n-\nu, a_1, \ldots ,a_p\\b_1, \ldots , b_q\end{array}\!;y\right],
\ee
where $p$ and $q$ are non-negative integers and ${}_pF_q$ denotes the generalized hypergeometric function.
In this, if we take $p=0$, $q=1$, $b_1=1+\nu+j$ and $y=-1$, then
\bee\label{e2}
e^{-x} \sum_{n=0}^\infty \frac{x^n\,L_n^{(\nu)}(x)}{(1+\nu+j)_n}=\sum_{n=0}^\infty \frac{(-x)^n}{n!}\,{}_2F_1\left[\!\!\begin{array}{c}-n, -n-\nu\\1+\nu+j\end{array}\!; -1\right].
\ee
The ${}_2F_1$ series on the right-hand side of (\ref{e2}) can be evaluated with the help of the generalized Kummer summation theorem \cite{RR}
\[{}_2F_1\left[\!\!\begin{array}{c}a, b\\1+a-b+j\end{array}\!;-1\right]\hspace{8cm}\]
\bee\label{e3}
=\frac{2^{-a}\g(\fs) \g(b-j) \g(1+a-b+j)}{\g(b) \g(\fs a-b+\fs j+\fs) \g(\fs a-b+\fs j+1)}\sum_{r=0}^j (-1)^r \left(\!\!\begin{array}{c}j\\r\end{array}\!\!\right)\,\frac{\g(\fs a-b+\fs j+\fs r+\fs)}{\g(\fs a-\fs j+\fs r+\fs)}
\ee
for $j=0, 1, 2, \ldots\ $. 
After some straightforward simplification, we obtain
\[S(\nu, j)\equiv e^{-x}\sum_{n=0}^\infty \frac{x^n\,L_n^{(\nu)}(x)}{(1+\nu+j)_n}\hspace{8cm}\]
\[
=\frac{(-1)^j 2^{2\nu+j} \g(1+\nu)}{\g(1+2\nu+j)} \sum_{r=0}^j (-1)^r \left(\!\!\begin{array}{c}j\\r\end{array}\!\!\right)\left\{\frac{\g(\nu+\fs j+\fs r+\fs)}{\g(\fs-\fs j+\fs r)}\right.\]
\[\hspace{2.5cm}\times {}_4F_5\left[\!\!\begin{array}{c}\vspace{0.1cm}

\fs+\fs\nu, 1+\fs\nu, \fs+\nu+\fs j+\fs r, \fs+\fs j-\fs r\\ \fs, \fs+\fs\nu+\fs j, 1+\fs\nu+\fs j, \fs+\nu+\fs j, 1+\nu+\fs j\end{array}\!;-x^2\right]\]
\[-\frac{4x(1+\nu)}{(1+\nu+j)(1+2\nu+j)}\,\frac{\g(\nu+\fs j+\fs r+1)}{\g(\fs r-\fs j)}\hspace{2cm}\]
\bee\label{e4}
\left.\hspace{3.2cm}\times {}_4F_5\left[\!\!\begin{array}{c}\vspace{0.1cm}

1+\fs\nu, \f{3}{2}+\fs\nu, 1+\nu+\fs j+\fs r, 1+\fs j-\fs r\\ \f{3}{2}, 1+\fs\nu+\fs j, \f{3}{2}+\fs\nu+\fs j, 1+\nu+\fs j, \f{3}{2}+\nu+\fs j\end{array}\!;-x^2\right]\right\}
\ee
for $j=0, 1, 2, \ldots\ $.

Again, in (\ref{e1}), if we take $p=0$, $q=1$, $b_1=1+\nu-j$ and $y=-1$, then 
\bee\label{e5}
e^{-x} \sum_{n=0}^\infty \frac{x^n\,L_n^{(\nu)}(x)}{(1+\nu-j)_n}=\sum_{n=0}^\infty \frac{(-x)^n}{n!}\,{}_2F_1\left[\!\!\begin{array}{c}-n, -n-\nu\\1+\nu-j\end{array}\!; -1\right].
\ee
The ${}_2F_1$ series on the right-hand side of (\ref{e5}) can be evaluated with the help of the known result \cite{RR}
\[{}_2F_1\left[\!\!\begin{array}{c}a, b\\1+a-b-j\end{array}\!;-1\right]\hspace{8cm}\]
\bee\label{e6}
=\frac{2^{-a}\g(\fs) \g(1+a-b-j)}{\g(\fs a-b-\fs j+\fs) \g(\fs a-b-\fs j+1)}\sum_{r=0}^j \left(\!\!\begin{array}{c}j\\r\end{array}\!\!\right)\,\frac{\g(\fs a-b-
\fs j+\fs r+\fs)}{\g(\fs a-\fs j+\fs r+\fs)}
\ee
for $j=0, 1, 2, \ldots$ and, after some simplification, we obtain
\[S(\nu, -j)\equiv e^{-x}\sum_{n=0}^\infty \frac{x^n\,L_n^{(\nu)}(x)}{(1+\nu-j)_n}\hspace{8cm}\]
\[
=\frac{2^{2\nu-j} \g(1+\nu-j)}{\g(1+2\nu-j)} \sum_{r=0}^j \left(\!\!\begin{array}{c}j\\r\end{array}\!\!\right)\left\{\frac{\g(\nu-\fs j+\fs r+\fs)}{\g(\fs-\fs j+\fs r)}\right.\]
\[\hspace{2.5cm}\times {}_2F_3\left[\!\!\begin{array}{c}\vspace{0.1cm}

\fs+\nu-\fs j+\fs r, \fs+\fs j-\fs r\\ \fs, \fs+\nu-\fs j, 1+\nu-\fs j\end{array}\!;-x^2\right]\]
\bee\label{e7}
\left.-\frac{4x}{(1+2\nu-j)}\,\frac{\g(\nu-\fs j+\fs r+1)}{\g(\fs r-\fs j)}
\times {}_2F_3\left[\!\!\begin{array}{c}\vspace{0.1cm}

1+\nu-\fs j+\fs r, 1+\fs j-\fs r\\ \f{3}{2}, 1+\nu-\fs j, \f{3}{2}+\nu-\fs j\end{array}\!;-x^2\right]\right\}
\ee
for $j=0, 1, 2, \ldots\ $.
\vspace{0.6cm}

\begin{center}
{\bf 3. \  THE SERIES $S(-\nu, \pm j)$}
\end{center}
\setcounter{section}{3}
\setcounter{equation}{0}
\renewcommand{\theequation}{\arabic{section}.\arabic{equation}}
Further, if we take $p=0$, $q=1$, $b_1=1-\nu+j$ and $y=-1$ in (\ref{e1}), we find
\bee\label{e8}
e^{-x} \sum_{n=0}^\infty \frac{x^n\,L_n^{(\nu)}(x)}{(1-\nu+j)_n}=\sum_{n=0}^\infty \frac{(-x)^n}{n!}\,{}_2F_1\left[\!\!\begin{array}{c}-n, -n-\nu\\1-\nu+j\end{array}\!; -1\right].
\ee
The ${}_2F_1$ series on the right-hand side of (\ref{e8}) can be evaluated by (\ref{e3})
to produce the result after some simplification
\[S(-\nu, j)\equiv e^{-x}\sum_{n=0}^\infty \frac{x^n\,L_n^{(\nu)}(x)}{(1-\nu+j)_n}\hspace{8cm}\]
\[
=\frac{(-2)^j}{j!} \sum_{r=0}^j (-1)^r\left(\!\!\begin{array}{c}j\\r\end{array}\!\!\right)\left\{\frac{\g(-\fs\nu+\fs j+\fs r+\fs)}{\g(-\fs\nu-\fs j+\fs r+\fs)}\hspace{1.5cm}\right.\]
\[\hspace{2.5cm}\times {}_3F_4\left[\!\!\begin{array}{c}\vspace{0.1cm}

1, \fs+\fs\nu+\fs j-\fs r, \fs-\fs \nu+\fs j+\fs r\\ \fs+\fs j, 1+\fs j, \fs-\fs\nu+\fs j, 1-\fs\nu+\fs j\end{array}\!;-x^2\right]\]
\[
-\frac{4x}{(j+1)(1-\nu+j)}\,\frac{\g(-\fs\nu+\fs j+\fs r+1)}{\g(-\fs\nu-\fs j+\fs r)}\hspace{2cm}\]
\bee\label{e9}
\hspace{2.5cm}\left.\times {}_3F_4\left[\!\!\begin{array}{c}\vspace{0.1cm}

1, 1+\fs\nu+\fs j-\fs r, 1-\fs\nu+\fs j+\fs r\\1+ \fs j, \f{3}{2}+\fs j, 1-\fs\nu+\fs j, \f{3}{2}-\fs\nu+\fs j\end{array}\!;-x^2\right]\right\}
\ee
for $j=0, 1, 2, \ldots\ $. 

Finally, if we take $p=0$, $q=1$, $b_1=1-\nu-j$ and $y=-1$ in (\ref{e1}), we find
\bee\label{e10}
 e^{-x} \sum_{n=0}^\infty \frac{x^n\,L_n^{(\nu)}(x)}{(1-\nu-j)_n}=\sum_{n=0}^\infty \frac{(-x)^n}{n!}\,{}_2F_1\left[\!\!\begin{array}{c}-n, -n-\nu\\1-\nu-j\end{array}\!; -1\right].
\ee
The ${}_2F_1$ series on the right-hand side of (\ref{e10}) can be evaluated by (\ref{e6})
to produce the result after some simplification
\[S(-\nu, -j)\equiv e^{-x}\sum_{n=0}^\infty \frac{x^n\,L_n^{(\nu)}(x)}{(1-\nu-j)_n}\hspace{8cm}\]
\[=2^{-j}\sum_{r=0}^j \left(\!\!\begin{array}{c}j\\r\end{array}\!\!\right)\left\{
{}_2F_3\left[\!\!\begin{array}{c}\vspace{0.1cm}

\fs-\fs\nu-\fs j+\fs r, \fs+\fs\nu+\fs j-\fs r\\ \fs, \fs-\fs\nu-\fs j, 1-\fs\nu-\fs j\end{array}\!;-x^2\right]\right.\]
\bee\label{e11}
\left.-\frac{2x(\nu+j-r)}{\nu+j-1}
{}_2F_3\left[\!\!\begin{array}{c}\vspace{0.1cm}

 1-\fs\nu-\fs j+\fs r, 1+\fs\nu+\fs j-\fs r\\ \f{3}{2}, 1-\fs\nu-\fs j, \f{3}{2}-\fs\nu-\fs j\end{array}\!;-x^2\right]\right\}
\ee
for $j=0, 1, 2, \ldots\ $. 
\vspace{0.6cm}

\begin{center}
{\bf 4. \  CONCLUDING REMARKS}
\end{center}
\setcounter{section}{4}
\setcounter{equation}{0}
\renewcommand{\theequation}{\arabic{section}.\arabic{equation}}
To conclude we make a brief comparison of the results (\ref{e4}), (\ref{e7}), (\ref{e9}) and (\ref{e11}) with those obtained in \cite{B}. The summations $S(\nu, \pm j)$ derived by Brychkov were expressed respectively in terms of finite sums of ${}_2F_3(-x^2)$ functions and Bessel functions of the first kind. The summations $S(-\nu, \pm j)$ were expressed respectively in terms of finite sums of four ${}_4F_3(-x^2)$ functions and four ${}_6F_7(-x^2)$ functions,
including the Jacobi polynomials of zero argument. Our expressions in (\ref{e9}) and (\ref{e11}) involve simpler finite sums of two ${}_3F_4(-x^2)$ and two ${}_2F_3(-x^2)$ functions, respectively. 

Finally, we mention that the summations $S(\pm\nu, \pm j)$ have been verified numerically with the help of {\it Mathematica}.

\vspace{0.8cm}

\noindent{\bf ACKNOWLEDGEMENT:}\ \ \ One of the authors (YSK) acknowledges the support of the Wonkwang University Research Fund (2013).

\end{document}